\documentclass[11pt,a4paper]{article}

\usepackage{amsmath}
\usepackage{graphicx,color}
\usepackage{graphpap,color}
\usepackage[mathscr]{eucal}
\usepackage{graphicx}
\usepackage{makeidx}
\usepackage{amsfonts}
\usepackage{latexsym}
\usepackage{amssymb}
\usepackage{float}

\newtheorem{t1}{Theorem}[section]

\newtheorem{l1}{Lemma}[section]
\newtheorem{conj}{Conjecture}[section]

\numberwithin{equation}{section}

\numberwithin{figure}{section}

\begin{document}
\begin{center}
{\Large \textbf{The asymptotic distribution of the number of 3-star factors in random \emph{d}-regular graphs}}\\
\bigskip
\vspace{1cm} Lyra Yulianti\footnote{This research partially
supported by the National Higher Education Grant No.
437/SP2H/PP/D2PM/V/2009}\\
\vspace{0.3 cm} {\small Department of Mathematics,\\Faculty of
Mathematics and Natural Sciences, Andalas University,\\
Kampus UNAND Limau Manis Padang 25136, Indonesia}\\
{\small lyra@fmipa.unand.ac.id}
\end{center}
\vspace{0.3 cm}
\begin{abstract}
The Small Subgraph Conditioning Method has been used to study the
almost sure existence and the asymptotic distribution of the number
of regular spanning subgraphs of various types in random
\emph{d}-regular graphs. In this paper we use the method to
determine the asymptotic distribution of the number of 3-star
factors in random \emph{d}-regular graphs for $d \geq 4$. \\
\textbf{Keywords} : Random regular graphs, Small Subgraph
Conditioning
Method\\
\textbf{AMS SUBJECT} : 05C80
\end{abstract}

\section{Introduction}
It their remarkable papers (see \cite{Wor92} and \cite{Wor94})
Robinson and Wormald showed that for $d \geq 3$ and $dn$ even, a
random \emph{d}-regular graph contains a Hamilton cycle with
probability tends to 1 as the number \emph{n} of vertices tends to
infinity. They used the \emph{Small Subgraph Conditioning Method}
(see \cite{Jan00} or \cite{Wor97} for details) to prove the
existence (with high probability) of perfect matching in such graphs
when $n$ is even.

The method has been used to determine the existence with high
probabili-ty of, and the asymptotic distribution of, the number of
\emph{k}-regular spanning subgraphs (for $k=1, 2$) and the number of
long cycles in random \emph{d}-regular graphs (see \cite{Gar99},
\cite{Jan95} and \cite{Rob96}).

A \emph{star} is a tree with at most one vertex whose degree is
greater than 1. A \emph{k-star} is a star with \emph{k} leaves. A
\emph{k-star factor} in a graph is a spanning subgraph whose
components are \emph{k}-stars.

We use notations \textbf{P}(probability), \textbf{E}(expectation)and
\textbf{Var}(variance). We say that an event $Y_n$ occurs
\emph{a.a.s (asymptotically almost surely)} if \[\lim_{n \rightarrow
\infty}~\textbf{P}Y_n = 1.\]

In~\cite{A06} Assiyatun and Wormald have used the method to investigate
the a.a.s of 3-star factor in random $d$-regular graphs. This is the first time the method applied to
non-regular subgraphs in such graphs.

Assiyatun and Wormald~\cite{A06} started proving the existence a.a.s
of a 3-star factor in random \emph{d}-regular graphs for $d \geq 4$
by showing the existence a.a.s of a 3-star factor in random
4-regular graphs. Then using the contiguity of models of random
regular graphs (see \cite{Wor99}), they obtained the existence a.a.s
of a 3-star factor in random \emph{d}-regular graphs for fixed $d
\geq 4$ as desired.

As a completion to the result in~\cite{A06}, in this paper we use
the method to determine the asymptotic distribution of the number of
3-star factors in random \emph{d}-regular graphs for $d \geq 4$.
However, due to the complexity of some part of the computation we
are only able to obtain the asymptotic distribution for $4\leq d
\leq 10$. Nevertheless, in most part of the computation we obtain
the result for general $d \geq 4$. The main result obtained in this
paper is presented in the following theorem.

Let $\mathcal{G}_{n,d}$ be a probability space contains of
\emph{d}-regular graphs with \emph{n} vertices. In asymptotic
statements about properties of $\mathcal{G}_{n,d}$, we restrict
\emph{n} to even integers when \emph{d} is odd.
\begin{t1}\label{main}
Restrict $n$ to $0$ mod $4$ and $4\leq d \leq 10$. Then $G\in
\mathcal{G}_{n, d}$ a.a.s has a 3-star factor. Furthermore, letting
$Y_d$ denote the number of 3-star factors in $G \in
\mathcal{G}_{n,d}$,
\[\frac{Y_d}{\textbf{E}Y_d}\rightarrow W =\prod_{k=3}^\infty(1+\delta_k)^{Z_k}
 e^{-\lambda_k \delta_k}\mbox{ for } n \rightarrow \infty\]
where $Z_k$ are independent Poisson variables with $\textbf{E}Z_k = \lambda_k$ for $k\geq 3$ and
\begin{eqnarray*}
\lambda_k & = & \frac{(d-1)^k}{2k}, \\
\delta_k & = & \left(\frac{-3(d-2)+\sqrt{-15d^2+24d}}{4(d-1)(d-3/2)}\right)^k
+\left(\frac{-3(d-2)-\sqrt{-15d^2+24d}}{4(d-1)(d-3/2)}\right)^k.
\end{eqnarray*}
\end{t1}

As in~\cite{Rob96} we will first work on the pairing model which was
first introduced by Bollob\'as (see \cite{Bol85}). This model can be
described as follows. Let $V = \bigcup_{i = 1}^n V_i$ be a fixed set
of \emph{dn} points, where $|V_i| = d$ for every \emph{i}. A
\emph{pairing} is defined as a perfect matching of points of $V$
into $dn/2$ pairs. A pairing \emph{P} corresponds to a random
\emph{d}-regular pseudograph \emph{G(P)} in which every $V_i$ is
regarded as a vertex and each pair is an edge. We use
$\mathcal{P}_{n,d}$ to denote the probability space of all pairings.

As shown in \cite{Bol85}, the probability that the pseudograph has
no loops or multiple edges (i.e simple graph) for a fixed \emph{d},
is asymptotically bounded below by a positive constant. Moreover,
each simple graph arises with the same probability as $G(P)$ for $P
\in \mathcal{P}_{n,d}$. Hence using the following property we obtain
the desired result in $\mathcal{G}_{n,d}$.
\begin{l1}\label{pairing}
A property of graphs that holds a.a.s for random pseudographs
arising from $\mathcal{P}_{n,d}$ will also hold a.a.s for
$\mathcal{G}_{n,d}$.
\end{l1}

Given two sequences $a_n$ and $b_n$, we denote $a_n \sim b_n$ if
${\large\frac{a_n}{b_n}} \to 1$ for $n \to \infty$. We denote
falling factorial $n(n-1)\cdots(n-m+1)$ by $[n]_m$ and Stirling's
formula by
\[n! \sim \sqrt{2\pi n}~\left(\frac{n}{e}\right)^n \mbox{  for  } n \to \infty.\]

\section{The variance of the number of star factors}
Throughout this paper we define
\[N(2m)=\displaystyle{\frac{(2m)!}{m!~2^m}}\]
as the number of perfect matchings of $2m$ points.

Note that counting subgraphs of the pseudograph coming from
$\mathcal{P}_{n,d}$ is equivalent to counting the corresponding sets
of pairs in the pairing. For that purpose, parallel edges are
distinguishable from each other (especially as they come from
distinct pairs in the pairing).

Let $n\equiv 0$ mod 4 and define $Y_d^*$ as the number of 3-star
factors in $G(P)$ coming from $\mathcal{P}_{n,d}$ for $d\geq 4$. We
have the following theorem.
\begin{t1}$\cite{A06}$\label{expectation}
\begin{equation}\label{ekspektasi}
\textbf{E}Y_d^* \sim 2\left(d (d-3/2)^{d/2-3/4}\left(\frac{2}{d^d}\right)^{1/2}
\left(\frac{(d-1)(d-2)}{3!}\right)^{1/4}\right)^n.
\end{equation}
\end{t1}
\textbf{Proof.} See \cite{A06} for details. $\square$\\

Using the method in \cite{A06} (see \cite{Fri96} and \cite{Rob96}
for similar argument) we obtain the following theorem.
\begin{t1}\label{variance}$\cite{L06}$
Restrict $n$ to $0$ mod $4$ and define $Y_d^*$ as the number of
3-star factors in $G(P)$ coming from $\mathcal{P}_{n,d}$ for $4 \leq
d \leq 10$. Then
\begin{equation}
\textbf{Var}Y_d^*\sim\left(\frac{2(d-1)^{1/2}(d-3/2)^2}{(d-3)(4d^3-13d^2+36d-36)^{1/2}}-1\right)(\textbf{E}Y_d^*)^2.
\end{equation}
\end{t1}
\textbf{Proof.} We count the ways to lay down an ordered pair of
3-star factors in $P \in \mathcal{P}_{n,d}$. In general, a set of
pairs in $P$ inducing a subgraph of a given type will be called by
the same name in the pairing.

Let $S_i$ be a 3-star factor of $P$ for $i=1,2$. Let $T= S_1 \cap
S_2$. Define $T=\bigcup_{j=1}^5~T_j$, where (see Figure 2.1)
\begin{itemize}
\item [(i.)] $T_1$ consists of $x_1$ 1-stars, $S_1$ and
$S_2$ have only one common leaf, \item [(ii.)] $T_2$ consists of
$x_2$ 1-stars, $S_1$ and $S_2$ have one common leaf and one common
center, \item [(iii.)] $T_3$ consists of $x_3$ 2-stars, $S_1$ and
$S_2$ have two common leaves, \item [(iv.)] $T_4$ consists of $x_4$
3-stars, $S_1 \simeq S_2$,
\item [(v.)] $T_5$ consists of $x_5$ \emph{0}-stars, $S_1$ and
$S_2$ have only one common center.
\end{itemize}

\begin{figure}[htbp]
        \centerline{\includegraphics[width=9cm]{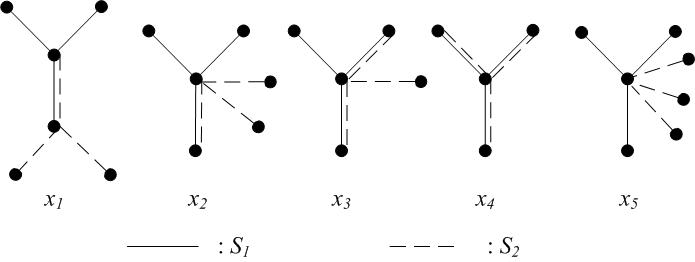}}
        \caption{Two intersecting 3-star factors}
\end{figure}

Given $S_1$, the number of possibilities of the intersection $T$ is
\[{n/4 \choose x_1}~3^{x_1}~{n/4- x_1 \choose x_2}~3^{x_2}~{n/4-x_1-x_2 \choose x_3}~3^{x_3}
~{n/4-x_1-x_2-x_3 \choose x_4}~1^{x_4}~~~~~~~~~~~~~~~~~~~~~~~~~~~\]
\begin{equation}\label{eq3}
\times {n/4-x_1-x_2-x_3-x_4 \choose x_5}~ =\frac{(n/4)!~3^{x_1+x_2+x_3}}{(n/4-x_1-x_2-x_3-x_4-x_5)!~x_1!x_2!x_3!x_4!x_5!}.
\end{equation}

There are $(\frac{n}{4}-x_1-x_2-x_3-x_4-x_5)$ 3-stars in both 3-star
factors that don't share any edge. We call these edge-disjoint
3-stars \emph{isolated 3-stars}.

We have to complete $S_2$ by creating the isolated 3-stars and
completing $T_1, T_2, T_3$ and $T_5$ into 3-stars.

The centers of the isolated 3-stars in $S_2$ can not be chosen from
the vertex set of $T$ nor the center of the isolated 3-stars in
$S_1$. There are
$(2x_1+2x_2+3x_3+4x_4+x_5)+(n/4-x_1-x_2-x_3-x_4-x_5)=(n/4+x_1+x_2+2x_3+3x_4)$
such vertices. Thus the number of ways to choose these centers is
\begin{equation}
{{3n/4-x_1-x_2-2x_3-3x_4} \choose {n/4-x_1-x_2-x_3-x_4-x_5}}
=\frac{(3n/4-x_1-x_2-2x_3-3x_4)!}{(n/4-x_1-x_2-x_3-x_4-x_5)!~(n/2-x_3-2x_4+x_5)!}.
\end{equation}

The number of ways to choose the leaves of the isolated 3-stars in
$S_2$ is
\[\prod_{k=0}^{n/4-x_1-x_2-x_3-x_4-x_5-1}{{3n/4-x_1-x_2-2x_3-3x_4-3k} \choose {3}}\]
\begin{equation}
=\frac{(3n/4-x_1-x_2-2x_3-3x_4)!}{(3!)^{n/4-x_1-x_2-x_3-x_4-x_5}~(2x_1+2x_2+x_3+3x_5)!}.
\end{equation}

The number of ways to choose the leaves of $T_5$ is
\begin{equation}
\prod_{k=0}^{x_5-1}{{2x_1+2x_2+x_3+3x_5-3k} \choose {3}}=\frac{(2x_1+2x_2+x_3+3x_5)!}{(3!)^{x_5}~(2x_1+2x_2+x_3)!}.
\end{equation}

The number of ways to choose the leaves for the completion of $T_3$ is
\begin{equation}
\prod_{k=0}^{x_3-1}{{2x_1+2x_2+x_3-k} \choose {1}}=\frac{(2x_1+2x_2+x_3)!}{(2x_1+2x_2)!}.
\end{equation}

For completing $T_2$ and $T_1$, the number of ways to choose the leaves are given consecutively
\begin{equation}
\prod_{k=0}^{x_2-1}{{2x_1+2x_2-2k} \choose {2}}=\frac{(2x_1+2x_2)!}{(2!)^{x_2}~(2x_1)!}
\end{equation}
and
\begin{equation}
\prod_{k=0}^{x_1-1}{{2x_1-2k} \choose {2}}=\frac{(2x_1)!}{(2!)^{x_1}}.
\end{equation}

So far we have determined the graph corresponding to $S_2$ but not
chosen the pairs of points corresponding to its edges. The number
of choices for these points is
\[(d-1)^{3n/4-x_2-2x_3-3x_4}~(d-2)^{n/4-x_2-x_3-x_4-x_5}~(d-3)^{n/2-2x_1-x_2-x_3-2x_4-x_5}\]
\begin{equation}
 \times~(d-4)^{x_2+x_5}~(d-5)^{x_5}.~~~~~~~~~~~~~~~~~~~~~~~~~~~~~~~~~~~~~~~~~~~~~~~~~~~~~~~~~~~~~~~~~~~~~~~~
\end{equation}
\indent Having $S_1$ and $S_2$, we observe that in $T_1$ there are
$2x_1$ vertices of degree 3, while in $T_2$ there are $x_2$ vertices
of degree 1 and $x_2$ vertices of degree 5. In $T_3$ we have $2x_3$
vertices of degree 1 and $x_3$ vertices of degree 4, while in $T_4$
we have $3x_4$ vertices of degree 1 and $x_4$ vertices of degree 3
and in $T_5$ we have $x_5$ vertices of degree 6. For the isolated
3-stars in $S_1$ and $S_2$ we have
$2(n/4-x_2-x_3-x_4-x_5)=n/2-2x_2-2x_3-2x_4-2x_5$ vertices of degree
4. The remaining vertices (there are $n/2-x_3-2x_4+x_5$ vertices)
are of degree 2. Thus the number of free points in $V$ is
{\setlength\arraycolsep{1pt}
\begin{eqnarray*}
&&2(d-3)x_1+(d-1)x_2+(d-5)x_2+2(d-1)x_3+(d-4)x_3+3(d-1)x_4+(d-3)x_4\\
&&+~(d-6)x_5+(d-4)(n/2-2x_2-2x_3-2x_4-2x_5)+(d-2)(n/2-x_3-2x_4+x_5)\\
&&=(d-3)n+2x_1+2x_2+4x_3+6x_4.
\end{eqnarray*}
\noindent Therefore the number of ways to complete the pairing $P$ is
\begin{equation}\label{eq4}
N((d-3)n+2x_1+2x_2+4x_3+6x_4).
\end{equation}

Multiplying equations (\ref{eq3}) -- (\ref{eq4}) by the number of
ways to choose $S_1$ as in (\ref{ekspektasi}), then dividing by
$N(nd)$, we have {\setlength\arraycolsep{2pt}
\begin{eqnarray*}
\textbf{E}Y_d^*(Y_d^*-1)&=&\frac{n!\left(\frac{nd}{2}\right)!}{(nd)!}\left(\frac{2}{\sqrt3}\right)^n
\left(d(d-1)\left((d-2)(d-3)\right)^{1/2}\right)^n\\
&& \times~ \sum_{R_1} \left(\frac{(3n/4-x_1-x_2-2x_3-3x_4)!^2}
{(n/4-x_1-x_2-x_3-x_4-x_5)!^2~(n/2-x_3-2x_4+x_5)!}\right.\\
&& \times~\frac{3^{2x_1+2x_2+2x_3+x_4}~((d-3)n+2x_1+2x_2+4x_3+6x_4)!}{2^{x_1+x_2+x_3+2x_4}~x_1!x_2!x_3!x_4!x_5!
~\left(\left(\frac{d-3}{2}\right)n+x_1+x_2+2x_3+3x_4\right)!}\\
&& \times~\left.\frac{(d-4)^{x_2+x_5}~(d-5)^{x_5}}{
(d-1)^{x_2+2x_3+3x_4}~(d-2)^{x_2+x_3+x_4+x_5}~(d-3)^{2x_1+x_2+x_3+2x_4+x_5}}\right)
\end{eqnarray*}}
\noindent where
\[R_1 =\{(x_1,x_2,x_3,x_4,x_5)\mid x_1,x_2,x_3,x_4,x_5 \geq 0,x_1+x_2+x_3+x_4+x_5 \leq n/4\}.\]

Set
\begin{eqnarray*}
p=\frac{x_1}{n}, \qquad q=\frac{x_2}{n}, \qquad r=\frac{x_3}{n}, \qquad s=\frac{x_4}{n}, \qquad t=\frac{x_5}{n}.
\end{eqnarray*}

Then Stirling's formula gives
{\setlength\arraycolsep{2pt}
\begin{eqnarray}\label{eq5}
\textbf{E}Y_d^*(Y_d^*-1) & \sim & \frac{1}{8(n \pi)^{5/2}}\left(\sqrt{2}~(2d)^{1-d/2}(d-1)
\left(\frac{(d-2)(d-3)}{3!}\right)^{1/2}\right)^n \nonumber \\
&& \times \sum_{R_2} \alpha(p, q, r, s, t)(F(p, q, r, s, t))^n
\end{eqnarray}}
where
{\setlength\arraycolsep{2pt}
\begin{eqnarray*}
R_2&=&\{(p, q, r, s, t)\mid p, q, r, s, t \geq 0, p+q+r+s+t\leq 1/4\}\\
F(p, q, r, s, t)&=&
\frac{f(3/4-p-q-2r-3s)^2~f\left(\left(\frac{d-3}{2}\right)+p+q+2r+3s\right)}{f(1/4-p-q-r-s-t)^2~
f(1/2-r-2s+t)~f(p)f(q)f(r)f(s)f(t)}\\
&& \times~
\frac{3^{2p+2q+2r+s}~2^{d-3+p+q+3r+4s}~(d-4)^{q+t}~(d-5)^t}{(d-1)^{q+2r+3s}~(d-2)^{q+r+s+t}~(d-3)^{2p+q+r+2s+t}}
\end{eqnarray*}}
with $f(x)=x^x$ and
\[\alpha(p, q, r, s, t)=\left(\frac{2(3/4-p-q-2r-3s)^2}{(1/4-p-q-r-s-t)^2~(1/2-r-2s+t)~pqrst}\right)^{1/2}.\]

Since by convention $f(0) := 1$, it can be seen that $F$ is
continuous in $R_2$. Next we determine the main contribution of the
sum which comes from the maximum of $F$ in $R_2$. The following
three lemmas prove that the maximum of $F$ is attained at
\[\textbf{x}_{d-max}=\left(\frac{9}{16d}, \frac{9(d-3)(d-4)}{16d(d-1)(d-2)}, \frac{9(d-3)}{8d(d-1)(d-2)}, \frac{3}{8d(d-1)(d-2)}, \frac{(d-3)(d-4)(d-5)}{16d(d-1)(d-2)}\right)\]
with
\[F(\textbf{x}_{d-max})=\frac{(2d)^{1-d/2}~(2(d-3/2))^{d-3/2}}{((d-1)(d-3))^{1/2}}.\]
\begin{l1}\label{xmax}$\cite{L06}$
Let $F$ and $R_2$ be as in $(\ref{eq5})$. Then for $d \geq 4$,
\[\textbf{x}_{d-max}=\left(\frac{9}{16d}, \frac{9(d-3)(d-4)}{16d(d-1)(d-2)}, \frac{9(d-3)}{8d(d-1)(d-2)}, \frac{3}{8d(d-1)(d-2)}, \frac{(d-3)(d-4)(d-5)}{16d(d-1)(d-2)}\right)\]
is the local maximum point of $F$ in the interior of $R_2$ with
\[F(\textbf{x}_{d-max})=\frac{(2d)^{1-d/2}~(2(d-3/2))^{d-3/2}}{((d-1)(d-3))^{1/2}}.\]
Moreover for $4 \leq d \leq 10$, $\textbf{x}_{d-max}$ is the
global maximum point of $F$ in the interior of $R_2$.
\end{l1}
\textbf{Proof.} First we look for all critical points of $F$ in the
interior of $R_2$. We set the partial derivations of $\ln{F}$ with
respect to $p, q, r, s$ and $t$, equal to 0, resulting in five
equations {\setlength\arraycolsep{1.5pt}
\begin{eqnarray}
0&=&9(1/4-p-q-r-s-t)^2(d-3+2p+2q+4r+6s)\nonumber\\
&&-p(d-3)^2(3/4-p-q-2r-3s)^2,\label{eqa}\\
0&=&9(1/4-p-q-r-s-t)^2(d-3+2p+2q+4r+6s)(d-4)\nonumber\\
&&-q(d-1)(d-2)(d-3)(3/4-p-q-2r-3s)^2,\label{eqb}\\
0&=&18(1/4-p-q-r-s-t)^2(d-3+2p+2q+4r+6s)^2(1/2-r-2s+t)\nonumber\\
&&-r(d-1)^2(d-2)(d-3)(3/4-p-q-2r-3s)^4,\label{eqc}\\
0&=&6(1/4-p-q-r-s-t)^2(d-3+2p+2q+4r+6s)^3(1/2-r-2s+t)^2\nonumber\\
&&-s(d-1)^3(d-2)(d-3)^2(3/4-p-q-2r-3s)^6,\label{eqd}\\
0&=&(1/4-p-q-r-s-t)^2(d-4)(d-5)-t(d-2)(d-3)(1/2-r-2s+t).\nonumber\\
&&\label{eqe}
\end{eqnarray}}
After substituting (\ref{eqa}) to (\ref{eqb}) and (\ref{eqc}) to (\ref{eqd}) we have
\begin{equation}
p=\frac{(d-1)(d-2)q}{(d-3)(d-4)},\label{eqp}
\end{equation}
and
\begin{eqnarray}
t&=&\frac{1}{r\left(d-3+{\displaystyle \frac{2(d-1)(d-2)q}{(d-3)(d-4)}}+2q+4r+6s\right)}\nonumber \\
&& \times \left(3s(d-1)(d-3)\left(3/4-\frac{(d-1)(d-2)q}{(d-3)(d-4)}-q-2r-3s\right)^2\right. \nonumber \\
&& \left.-~(1/2-r-2s)r\left(d-3+\frac{2(d-1)(d-2)q}{(d-3)(d-4)}+2q+4r+6s\right)\right).\nonumber\\
&& \label{eqt}
\end{eqnarray}}
After substituting (\ref{eqp}) and (\ref{eqt}) to (\ref{eqa}),
(\ref{eqc}) and (\ref{eqe}) we have three homogenous equations
\[\frac{P_1(q, r, s)}{Q_1(q, r, s)}=0,~~\frac{P_2(q, r, s)}{Q_2(q, r, s)}=0,~~\frac{P_3(q, r, s)}{Q_3(q, r, s)}=0,\]
where
\begin{eqnarray*}
Q_1(q, r, s)&=&256(d-3)^3(d-4)^3\left((4d^2-20d+28)q+(4d^2-28d+48)r+(6d^2-42d+72)s\right.\\
&&\left.+~d^3-10d^2+33d-36\right)r^2,\\
Q_2(q, r, s)&=&8(d-3)(d-4)^2r~Q_1(q, r, s),\\
Q_3(q, r, s)&=&\frac{(Q_1(q, r, s))^2}{512~(d-3)^4~(d-4)^4~r^2}
\end{eqnarray*}
and degrees of $q$ in $P_1(q, r, s), P_2(q, r, s)$ and $P_3(q, r,
s)$ are 2, while degrees of $r$ and $s$ are $5$. It is sufficient to
look at their nominator parts,
\begin{eqnarray}
P_1(q, r, s)&=&0,\label{res1}\\
P_2(q, r, s)&=&0,\label{res2}\\
P_3(q, r, s)&=&0.\label{res3}
\end{eqnarray}
After taking the resultant of (\ref{res1}) and (\ref{res2}) and of
(\ref{res1}) and (\ref{res3}) with respect to $q$, we have two
homogenous equations,
\begin{eqnarray}
U_1(r, s)~T_1(r, s)&=&0,\label{res4}\\
U_2(r, s)~T_2(r, s)&=&0,\label{res5}\
\end{eqnarray}
where
\begin{eqnarray*}
U_1(r, s)&=&-331776r^4s^2(d-2)^2(2d-3)^2(d^2-5d+7)^2(d-1)^4(d-4)^4(d-3)^6,\\
U_2(r, s)&=&\frac{(d-3)^4}{9r^2}~U_1(r, s)
\end{eqnarray*}
and degrees of $r$ and $s$ in $T_1(r, s)$ and $T_2(r, s)$ are $6$.
By taking the resultant of $T_1(r, s)$ and $T_2(r, s)$ with respect
to $r$, we have {\setlength\arraycolsep{1.5pt}
\begin{eqnarray}
0&=&-4497760410984972288s^{16}(d-1)^4(d-3)^{12}(4s-1)^4\nonumber\\
 & &(8d^3s-24sd^2+16sd-3)~V(s).\label{resakhir}
\end{eqnarray}}
where $V(s)$ is a polynomial of degree $15$ (see \cite{L06} for
details). It is easy to show that one of the feasible solution for
(\ref{resakhir}) is
\begin{equation}\label{smax}
s^*=\frac{3}{8d(d-1)(d-2)}.
\end{equation}
By substituting (\ref{smax}) to (\ref{res4}) and (\ref{res5}) we have
\begin{eqnarray}
W(r, s^*)~N_1(r, s^*)&=&0,\label{res4a}\\
W(r, s^*)~N_2(r, s^*)&=&0.\label{res5a}
\end{eqnarray}
where
\[W(r, s^*)=8rd^3-24rd^2+16rd+27-9d\]
and degrees of $r$ in $N_1(r, s^*)$ and $N_2(r, s^*)$ are 6. From
$W(r, s^*) = 0$ we have a feasible solution for $r$
\begin{equation}\label{rmax}
r^*=\frac{9(d-3)}{8d(d-1)(d-2)}.
\end{equation}
By substituting (\ref{rmax}) to (\ref{res1}) -- (\ref{res3}) we have three new equations
\begin{eqnarray}
Z(q, r^*, s^*)~S_1(q, r^*, s^*)&=&0,\label{res1a}\\
Z(q, r^*, s^*)~S_2(q, r^*, s^*)&=&0,\label{res2a}\\
Z(q, r^*, s^*)~S_3(q, r^*, s^*)&=&0.\label{res3a}
\end{eqnarray}
where \[Z(q, r^*, s^*)=16qd^3-9d^2-48qd^2+32qd+63d-108\] and degrees
of $q$ in $S_1(q, r^*, s^*)$, $S_2(q, r^*, s^*)$ and $S_3(q, r^*,
s^*)$ are 2. From $Z(q, r^*, s^*)=0$ we have a feasible solution for
$q$
\begin{equation}\label{qmax}
q^*=\frac{9(d-3)(d-4)}{16d(d-1)(d-2)}.
\end{equation}
By substituting (\ref{qmax}) to (\ref{eqp}) and (\ref{eqt}) we have
\begin{eqnarray*}
p^*&=&\frac{9}{16d},\\
t^*&=&\frac{(d-3)(d-4)(d-5)}{16d(d-1)(d-2)}.
\end{eqnarray*}
Then we obtain
\[\textbf{x}_{d-max}=\left(\frac{9}{16d}, \frac{9(d-3)(d-4)}{16d(d-1)(d-2)},
 \frac{9(d-3)}{8d(d-1)(d-2)},  \frac{3}{8d(d-1)(d-2)}, \frac{(d-3)(d-4)(d-5)}{16d(d-1)(d-2)}\right).\]
\textbf{Note:} See \cite{L06} for details of $P_1(q, r, s),~P_2(q, r, s),~P_3(q, r, s),~T_1(r, s),~T_2(r, s),~N_1(r, s^*)$,
 $N_2(r, s^*)$, $S_1(q, r^*, s^*)$, $S_2(q, r^*, s^*)$, $S_3(q, r^*, s^*)$ and $V(s)$.\\

For $d \geq 4$, the Hessian matrix of $F$ at $\textbf{x}_{d-max}$
is negative definite. Then $\textbf{x}_{d-max}$ is the local
maximum point of $F$ in the interior $R_2$ with
\[F(\textbf{x}_{d-max})=\frac{(2d)^{1-d/2}~(2(d-3/2))^{d-3/2}}{((d-1)(d-3))^{1/2}}.\]

To prove that $\textbf{x}_{d-max}$ is the global maximum point in
the interior of $R_2$ for $4\leq d\leq 10$, we use the following
procedure,
\begin{enumerate}
\item [(1)] For each $d$, determine the feasible solutions of
(\ref{resakhir}). If there are $k$ feasible solutions then denote
them by $s_{d1}, \ldots, s_{dk}$, \item [(2)] For $s=s_{di}$, $1
\leq i \leq k$, determine the feasible solutions of (\ref{res4})
-- (\ref{res5}). If there are $l_i$ feasible solutions for each
$i$, then denote them by $r_{di1}, \ldots, r_{dil_i}$, \item [(3)]
For $s=s_{di}$ and $r=r_{dij}$, $1 \leq i \leq k$, $1 \leq j \leq
l_i$, determine the feasible solutions of (\ref{res1}) --
(\ref{res3}). If there are $m_{ij}$ feasible solutions for each
$i$ and $j$, then denote them by $q_{dij1}, \ldots,
q_{dijm_{ij}}$, \item [(4)] For $s=s_{di}$, $r=r_{dij}$, and
$q=q_{dijf}$, $1 \leq i \leq k$, $1 \leq j \leq l_i$, $1 \leq f
\leq m_{ij}$ determine $p_{dijf}$ and $t_{dijf}$ in (\ref{eqp})
and (\ref{eqt}), \item [(5)] Define $\textbf{x}_{dijf}=(p_{dijf},
q_{dijf}, r_{dij}, s_{di}, t_{dijf})$ as another feasible solution
for system (\ref{eqa}) -- (\ref{eqe}) in the interior of $R_2$,
\item [(6)] Determine $F(\textbf{x}_{d-max}) -
F(\textbf{x}_{dijf})$.
\end{enumerate}

For each $d$, $4\leq d\leq 10$, we have $k = 1, l_i = 1, m_{ij} =
1$. It means that we have only one other feasible solution for (\ref{eqa}) -- (\ref{eqe})
 in the interior of $R_2$. Because $F(\textbf{x}_{dijf})< F(\textbf{x}_{d-max})$
  for each $d$, $4 \leq d \leq 10$, then $\textbf{x}_{d-max}$ is the global maximum point of $F$ in the interior of $R_2$. $\square$\\

To study the behavior of $F$ on the boundary of $R_2$, we generalize
the approach used by Garmo in the proof of (\cite{Gar99}, Lemma 12).
First let $\textbf{x}=(x_1, x_2, \ldots, x_r)$ and
$\textbf{u}_i=(u_{1, i}, u_{2, i}, \ldots, u_{r, i})$ for $r\geq2$.
The $\ln$ function is defined on the set of non-negative real
numbers with, by convention, $0\times \ln{0}=0$.
\begin{l1}$\cite{A06}$ \label{le2}
Let $R: = R_2$ be a closed set in $\mathbb{R}^r$ and let $\delta R$
be the boundary of $R$. Assume that every point in $\delta R$ is the
endpoint of an interval in $R \backslash \delta R$. Let
$f_i(\textbf{x})=b_i+\textbf{u}_i \textbf{x}^T$ for $i=1,\cdots ,m$,
where $b_i$ and $u_i$ are constant, such that $f_i(\textbf{x})>0$
for all $i$ and all $\textbf{x}\in R \backslash \delta R$. Define
$F$ to be a function on $R$ such that
\[F(\textbf{x})=g_0(\textbf{x})+\sum_{i=1}^{m} a_ig_i(\textbf{x})=
g_0(\textbf{x})+\sum_{i=1}^{m} a_if_i(\textbf{x})\ln
{f_i(\textbf{x})}.\] with $a_i < 0$ for $i \leq m_0 \leq m$. Suppose
that for every $\textbf{x} \in R$, the directional derivative of
$g_0$ at $\textbf{x}$ in any direction is bounded. Let
$\textbf{x}_0\in \delta R$ such that $f_i(\textbf{x}_0)=0$ for at
least one $i\leq m_0$ and $f_i(\textbf{x}_0)>0$ for all $m_0 < i
\leq m$. Then $\textbf{x}_0$ is not a local maximum of $F$ on $R$.
\end{l1}
\textbf{Proof.} See \cite{A06} or \cite{L06} for details.\\
\begin{l1}$\cite{L06}$
Let $F$ and $R:=R_2$ be as in $(\ref{eq5})$. Then the maximum of $F$ does not occur in $\delta R$.
\end{l1}
\textbf{Proof.} We define $\textbf{x}=(p, q, r, s, t)$ and
$\textbf{v}=(v_1, v_2, v_3, v_4, v_5)$. Following the notation in
Lemma \ref{le2} we write
\[\ln{F(\textbf{x})}=g_0(\textbf{x})+\sum_{i=1}^9 a_i f_i(\textbf{x})\ln{f_i(\textbf{x})}\]
where {\setlength\arraycolsep{2pt}
\begin{eqnarray*}
g_0(\textbf{x})&=&(d-3+2p+2q+4r+6s)\ln{(d-3+2p+2q+4r+6s)}+(2p+2q+2r+s)\ln{(3)}\\
&& +~(q+t)\ln{(d-4)}+(t)\ln{(d-5)}-(p+q+r+2s)\ln{(2)}-(q+2r+3s)\ln{(d-1)}\\
&& -~(q+r+s+t)\ln{(d-2)}-(2p+q+r+2s+t)\ln{(d-3)},
\end{eqnarray*}}
$a_1=a_2=a_3=a_4=a_5=a_6=a_8=-1$, $a_7=-2$, $a_9=2$, $b_1=b_2=b_3=b_4=b_5=0$,
 $b_6=\frac{d-3}{2}$, $b_7=1/4$, $b_8=1/2$, $b_9=3/4$ and
\begin{eqnarray*}
\textbf{u}_1&=&(1, 0, 0, 0, 0),\\
\textbf{u}_2&=&(0, 1, 0, 0, 0),\\
\textbf{u}_3&=&(0, 0, 1, 0, 0),\\
\textbf{u}_4&=&(0, 0, 0, 1, 0),\\
\textbf{u}_5&=&(0, 0, 0, 0, 1),\\
\textbf{u}_6&=&(1, 1, 2, 3, 0),\\
\textbf{u}_7&=&(-1, -1, -1, -1, -1),\\
\textbf{u}_8&=&(0, 0, -1, -2, 1),\\
\textbf{u}_9&=&(-1, -1, -2, -3, 0).
\end{eqnarray*}
For $g_0$ we have {\setlength\arraycolsep{1pt}
\begin{eqnarray*}
\left. \frac{\partial}{\partial c}~g_0(\textbf{x}+c\textbf{v})\right|_{c=0}
&=&(2v_1+2v_2+4v_3+6v_4)\ln{(d-3+2p+2q+4r+6s)}\\
&& +~(2v_1+2v_2+4v_3+6v_4)+(2v_1+2v_2+2v_3+v_4)\ln{(3)}\\
&& +~(v_2+v_5)\ln{(d-4)}+(v_5)\ln{(d-5)}\\
&& -~(v_1+v_2+v_3+2v_4)\ln{(2)}-(v_2+2v_3+3v_4)\ln{(d-1)}\\
&& -~(v_2+v_3+v_4+v_5)\ln{(d-2)}-(2v_1+v_2+v_3+2v_4+v_5)\ln{(d-3)}
\end{eqnarray*}}
which is bounded for all $\textbf{x}\in R$. Having $a_9>0$, then
from Lemma \ref{le2} we only need to consider the solution of the
following system
\begin{eqnarray*}
f_i(\textbf{x})& \geq & 0,~i=1, \ldots, 8\\
f_9(\textbf{x})& = &0
\end{eqnarray*}
which is equivalent to
\begin{eqnarray*}
p&\geq&0,\\
q&\geq&0,\\
r&\geq&0,\\
s&\geq&0,\\
t&\geq&0,\\
\frac{d-3}{2}+p+q+2r+3s &\geq&0,\\
1/4-p-q-r-s-t &\geq&0,\\
1/2-r-2s+t &\geq&0,\\
3/4-p-q-2r-3s&=&0.
\end{eqnarray*}
It is easy to show that the only solution for the system is
$\textbf{c}_1=(0, 0, 0, 1/4, 0)$. Consequently for $\textbf{x}_0 \in
\delta R \backslash \{\textbf{c}_1\}$, $\ln{F}$ and $x_0$ satisfies
the hypotheses of Lemma \ref{le2}. Hence $F$ does not have any
maximum on $\delta R \backslash \{\textbf{c}_1\}$. Moreover for
$d\geq 4$,
\[F(\textbf{c}_1)=\frac{3^{1/4}~(2(d-3/2))^{d/2-3/4}}{((d-1)^3(d-2)(d-3)^2)^{1/4}}\]
is strictly less than $F(\textbf{x}_{d-max})$. $\square$\\

In the following lemma we will show that the sum in (\ref{eq5}) can
be approximated within a small region around the maximum.
\begin{l1}\label{le3}$\cite{L06}$
Let $B=B(\textbf{x}_{d-max}, \delta)$ be a ball centered at
\[\textbf{x}_{d-max}=\left(\frac{9}{16d}, \frac{9(d-3)(d-4)}{16d(d-1)(d-2)}, \frac{9(d-3)}{8d(d-1)(d-2)},    \frac{3}{8d(d-1)(d-2)}, \frac{(d-3)(d-4)(d-5)}{16d(d-1)(d-2)}\right)\]
and diameter $\delta:=n^{-5/2}$. Then with $F$ and $R:=R_2$ as in $(\ref{eq5})$, we have
\[\sum_{R} \alpha (\textbf{x}) F^n(\textbf{x})\sim \sum_{B} \alpha (\textbf{x}) F^n(\textbf{x}).\]
\end{l1}
\textbf{Proof.} Write
\[\sum_{R} \alpha (\textbf{x}) F^n(\textbf{x})= \sum_{B} \alpha (\textbf{x}) F^n(\textbf{x})+\sum_{R \backslash B} \alpha (\textbf{x}) F^n(\textbf{x}).\]
It will be shown that
\[\sum_{R \backslash B} \alpha (\textbf{x}) F^n(\textbf{x})=o(\alpha(\textbf{x}_{d-max})F^n(\textbf{x}_{d-max})).\]
For $\textbf{x}\in B$, the Taylor expansion of $F$ at
$\textbf{x}_{d-max}$ is
\[F^n(\textbf{x})=F^n(\textbf{x}_{d-max})\times \sum_{B} \exp \left(-n\sum_{i=1}^{5}~\sum_{j=i}^5 c_{ij}s_is_j\right)\]
\noindent where
\[\left.\begin{array}{lll}
s_1=p-\displaystyle{\frac{9}{16d}},&
s_2=q-\displaystyle{\frac{9(d-3)(d-4)}{16d(d-1)(d-2)}}, &
s_3=r-\displaystyle{ \frac{9(d-3)}{8d(d-1)(d-2)}},\\
&&\\
s_4=s-\displaystyle{\frac{3}{8d(d-1)(d-2)}},&
s_5=t-\displaystyle{\frac{(d-3)(d-4)(d-5)}{16d(d-1)(d-2)}}&
\end{array}\right.\]
and $c_{ij}$ are coming from Hessian matrix, $H=(a_{ij})\in M_{5 x 5}$, of $F$ in the proof of Lemma \ref{xmax}
\[c_{ij}=\left\{ \begin{array}{ll}
  \frac{1}{2} a_{ij}, & \mbox{ if   } i=j\\
              a_{ij}, & \mbox{ if   } i\neq j
\end{array}\right.\]
For $\textbf{x}^*\in \delta B$, where $\delta B$ is the boundary of $B$, we note that
\[O(e^{n^{-1/5}})=o(1).\]
Then
\[\alpha(\textbf{x}^*)F^n(\textbf{x}^*)\sim \alpha(\textbf{x}_{d-max})F^n(\textbf{x}_{d-max})o(1)=
o(\alpha(\textbf{x}_{d-max})F^n(\textbf{x}_{d-max})).\] Since $F$
attains its maximum uniquely at $\textbf{x}_{d-max}$, then for
$\textbf{x} \in R \backslash B$
\[\alpha(\textbf{x})F^n(\textbf{x})=O\left(\max_{\textbf{x}^* \in \delta B} \alpha(\textbf{x}^*)F^n(\textbf{x}^*)\right).\]
Thus we have
\[\sum_{R\backslash B} \alpha(\textbf{x})F^n(\textbf{x})=o(\alpha(\textbf{x}_{d-max})F^n(\textbf{x}_{d-max})).~~~~~~~~~~~~~~~~~~~~~~~~\square\]

Now we determine $\sum_{B}$ $\alpha(\textbf{x})F^n(\textbf{x})$.
Since the summation concentrates near the maximum, each term
$\alpha(\textbf{x})$ can be taken as $\alpha(\textbf{x}_{d-max})$
with
\[\alpha(\textbf{x}_{d-max})=\frac{8192\sqrt{6}~d^3(d-1)^{3/2}(d-3/2)(d-2)^2}{243~(d-3)^{5/2}(d-4)(d-5)^{1/2}}.\]
Referring to the Taylor expansion of $F$ as in the proof of Lemma \ref{le3} we have
\[\sum_{B}\alpha(\textbf{x})F^n(\textbf{x})\sim \alpha(\textbf{x}_{d-max})F^n(\textbf{x}_{d-max})~\sum_{B}
\exp \left(-n\sum_{i=1}^{5}~\sum_{j=i}^5 c_{ij}s_is_j\right).\]
The summation is a Riemann integral for the five-tuple integral
\[n^{5/2}\int_{-n}^{n}\int_{-n}^{n}\int_{-n}^{n}\int_{-n}^{n}\int_{-n}^{n}~\exp \left(-n \sum_{i=1}^{5}~\sum_{j=i}^5~c_{ij}t_it_j\right)~dt_1\, dt_2\, dt_3\, dt_4\, dt_5\]
where
\[t_1=\displaystyle{\frac{\left(p-\displaystyle{\frac{9}{16d}}\right)n}{\sqrt{n}}},~~
t_2=\displaystyle{\frac{\left(q-\displaystyle{\frac{9(d-3)(d-4)}{16d(d-1)(d-2)}}\right)n}{\sqrt{n}}},~~
t_3=\displaystyle{\frac{\left(r-\displaystyle{\frac{9(d-3)}{8d(d-1)(d-2)}}\right)n}{\sqrt{n}}},\]
\[t_4=\displaystyle{\frac{\left(s-\displaystyle{\frac{3}{8d(d-1)(d-2)}}\right)n}{\sqrt{n}}},~~
t_5=\displaystyle{\frac{\left(t-\displaystyle{\frac{(d-3)(d-4)(d-5)}{16d(d-1)(d-2)}}\right)n}{\sqrt{n}}}.~~~~~~~~\]
As $n \rightarrow \infty$, the range of the integration can be extended to $\pm \infty$ without altering the main asymptotic term. Thus it is asymptotic to
\[n^{5/2}\int_{-\infty}^{\infty}\int_{-\infty}^{\infty}\int_{-\infty}^{\infty}\int_{-\infty}^{\infty}
\int_{-\infty}^{\infty}\exp
\left(-\left(\sum_{i=1}^{5}~\sum_{j=i}^5~c_{ij}~t_it_j\right)\right)~dt_1\,
dt_2\, dt_3\,dt_4\, dt_5.\]
The evaluation of the integral results in
\[\frac{162\sqrt{6}~\pi^{5/2}(d-3/2)(d-3)^{3/2}(d-4)(d-5)^{1/2}}{d^3(d-1)(d-2)^2(4d^3-13d^2+36d-36)^{1/2}}.\]\\
Then (\ref{eq5}) becomes
{\setlength\arraycolsep{1pt}
\begin{eqnarray*}
\textbf{E}Y_d^*(Y_d^*-1)&\sim& n^{5/2}\times \frac{1}{8(n\pi)^{5/2}}\left(\sqrt{2}~(2d)^{1-d/2}~(d-1)\left(\frac{(d-2)(d-3)}{3!}\right)^{1/2}\right)^n\\
&&\\
&&\times~\left(\frac{(2d)^{1-d/2}~(2(d-3/2))^{d-3/2}}{((d-1)(d-3))^{1/2}}\right)^n ~\frac{8192\sqrt{6}~d^3(d-1)^{3/2}(d-3/2)(d-2)^2}{243~(d-3)^{5/2}(d-4)(d-5)^{1/2}}\\
&&\\
&&\times~\frac{162\sqrt{6}~\pi^{5/2}(d-3/2)(d-3)^{3/2}(d-4)(d-5)^{1/2}}{d^3(d-1)(d-2)^2(4d^3-13d^2+36d-36)^{1/2}}
\end{eqnarray*}
\begin{eqnarray}\label{eq14}
\sim\frac{8(d-1)^{1/2}(d-3/2)^2}{(d-3)(4d^3-13d^2+36d-36)^{1/2}}\left(d^2(d-3/2)^{d-3/2}\left(\frac{2}{d^d}\right)\left(\frac{(d-1)(d-2)}{3!}\right)
^{1/2}\right)^n.&&\nonumber\\
\end{eqnarray}}
\indent As $\textbf{E}Y_d^* \rightarrow \infty$ for $d\geq 4$ we have $\textbf{E}Y_d^*(Y_d^*-1)\sim \textbf{E}Y_d^{*2}$. Thus from (\ref{eq14}) and (\ref{ekspektasi}) we have
\begin{equation}\label{eq15}
\textbf{E}Y_d^{*2} \sim
\frac{2(d-1)^{1/2}(d-3/2)^2}{(d-3)(4d^3-13d^2+36d-36)^{1/2}}~(\textbf{E}Y_d^*)^2.
\end{equation}\\
Since $\textbf{Var}Y_d^*=\textbf{E}Y_d^{*2}-(\textbf{E}Y_d^*)^2$,
the above equation gives the required result. $\square$

\section{Expectation conditioned on short cycle distribution}
\begin{l1}\label{le4}$\cite{L06}$
Let $n\equiv 0$ mod $4$ and $X_k$ be the number of cycles of
length $k$ in $G(P)$ for $P\in\mathcal{P}_{n,d}$. Then for any
finite sequences $j_1$, \ldots, $j_m$ of non-negative integers
\[\frac{\textbf{E}(Y_d^*[X_1]_{j_1} \ldots [X_m]_{j_m})}{\textbf{E}Y_d^*}\rightarrow \prod_{k=1}^m~(\lambda_k~
(1+\delta_k))^{j_k} \mbox{ for } n\rightarrow \infty\]
where
\begin{eqnarray*}
\lambda_k&=&\frac{(d-1)^k}{2k},\\
\delta_k&=&\left(\frac{-3(d-2)+\sqrt{-15d^2+24d}}{4(d-1)(d-3/2)}\right)^k
+\left(\frac{-3(d-2)-\sqrt{-15d^2+24d}}{4(d-1)(d-3/2)}\right)^k.
\end{eqnarray*}
\end{l1}
\textbf{Proof.} To prove the lemma we first establish
\begin{eqnarray}
\frac{\textbf{E}(Y_d^*X_k)}{\textbf{E}Y_d^*}& \sim &\frac{1}{2k} \left((d-1)^k+\left(\frac{-3(d-2)+\sqrt{-15d^2+24d}}{4(d-3/2)}\right)^k\right.\nonumber\\
&&\left.+~\left(\frac{-3(d-2)-\sqrt{-15d^2+24d}}{4(d-3/2)}\right)^k\right).\label{eq16}
\end{eqnarray}

The number of ways to choose a cycle of length $k$ in the pairing (with a distinguished point in a pair) is
\begin{equation}\label{eq17}
\frac{n!}{(n-k)!}~(d(d-1))^k.
\end{equation}

This induces an orientation and also a distinguished edge called a \emph{root edge} in the cycle.

Let $C$ denote the set of pairs that corresponds to an oriented and
rooted $k$-cycle. Define $S$ to be the set of pairs corresponding to
a 3-star factor. Fix $C$ and suppose $C\cap S$ consists of $s_0$
0-stars (By 0-star we mean isolated vertices) lying at the centers
of stars in the 3-star factor, $s_1$ 1-stars, $s_2$ 2-stars and
$s_3=k-s_0-2s_1-3s_2$ 0-stars lying at the leaves of the 3-star
factor.

The edges of $C$ can be classified into three types. The first type
is the edges not lying in the 3-star factor, we denote this as
\emph{$\textbf{0}$}. The second is the edges of 1-stars and the
first edges of 2-stars and the last is the second edges of 2-stars.
We denote them as \emph{$\textbf{1}$} and \emph{$\textbf{2}$}
respectively. If we walk along $C$ from the root edge, we obtain a
sequence
$S_0\in\{\emph{$\textbf{0}$},\emph{$\textbf{1}$},\emph{$\textbf{2}$}\}^k$.

Fix $C$ and $S_0$. The number of ways to choose the centers of the
remaining $(n/4-s_0-s_1-s_2)$ 3-stars, together with the points
used, is
\begin{equation}\label{eq17a}
{{n-k}\choose{n/4-s_0-s_1-s_2}}~(d(d-1)(d-2))^{n/4-s_0-s_1-s_2}.
\end{equation}

The number of ways to choose the points in the centers of
$(s_0+s_1+s_2)$ 3-stars is
\begin{equation}\label{eq18}
((d-2)(d-3)(d-4))^{s_0}((d-2)(d-3))^{s_1}(d-2)^{s_2}= (d-2)^{s_0+s_1+s_2}(d-3)^{s_0+s_1}(d-4)^{s_0}.
\end{equation}

The number of leaves remaining for the 3-star factor is
\[3(n/4-s_0-s_1-s_2)+3s_0+2s_1+s_2=3n/4-s_1-2s_2.\]
The number of ways to select the leaves from the remaining
$(n/4-s_0-s_1-s_2)$ 3-star is
\begin{equation}\label{eq19}
\prod_{k=0}^{n/4-s_0-s_1-s_2-1}
{{3n/4-s_1-2s_2-3k} \choose {3}}=\frac{(3n/4-s_1-2s_2)!}{(3!)^{n/4-s_0-s_1-s_2}(3s_0+2s_1+s_2)!}.~~~~~~~~
\end{equation}
The number of ways to choose the leaves from $s_1$ 1-stars
\begin{equation}
\prod_{k=0}^{s_1-1}
{{3s_0+2s_1+s_2-2k} \choose {2}}=\frac{(3s_0+2s_1+s_2)!}{(2!)^{s_1}(3s_0+s_2)!}.
\end{equation}
The number of ways to choose the leaves from $s_2$ \emph{2}-stars
\begin{equation}
\prod_{k=0}^{s_2-1}
{{3s_0+s_2-k} \choose {1}}=\frac{(3s_0+s_2)!}{(3s_0)!}.
\end{equation}
The number of ways to choose the leaves from $s_0$ 0-stars
\begin{equation}
\prod_{k=0}^{s_0-1}
{{3s_0-3k} \choose {3}}=\frac{(3s_0)!}{(3!)^{s_0}}.
\end{equation}
The number of ways to choose which vertex will be the center of $s_1$
\begin{equation}
2^{s_1}.
\end{equation}
The number of leaves of $S$ lying in the cycle is $k-s_0-2s_1-3s_2$. Note that every vertex in the cycle uses their
two points. Thus the number of ways to choose the points that represent the leaves is
\begin{equation}
(d-2)^{k-s_0-2s_1-3s_2}.
\end{equation}
There are $(n-k)-(n/4-s_0-s_1-s_2)=(3n/4-k+s_0+s_1+s_2)$ leaves
outside the cycle. Note that there are $d$ points in every vertex of
the cycle. Thus the number of ways choose the points that represent
the leaves is
\begin{equation}\label{eq20}
d^{3n/4-k+s_0+s_1+s_2}.
\end{equation}
By multiplying (\ref{eq19}) -- (\ref{eq20}) we have the number of ways to choose the leaves (including the points used)
\begin{equation}\label{eq21}
\frac{(3n/4-s_1-2s_2)!}{(3!)^{n/4-s_1-s_2}}~(d-2)^{k-s_0-2s_1-3s_2}~d^{3n/4-k+s_0+s_1+s_2}.
\end{equation}

The number of points for the leaves of 3-stars outside the cycles is
$(d-1)(3n/4-k+s_0+s_1+s_2)$, while for the centers is
$(d-3)(n/4-s_0-s_1-s_2)$. The number of points for the centers of
3-stars inside the cycle is  $(d-3)k-2s_0+2s_2$. Thus the number of
free points in the pairing is {\setlength\arraycolsep{1pt}
\begin{eqnarray*}
&&(d-1)(3n/4-k+s_0+s_1+s_2)+(d-3)(n/4-s_0-s_1-s_2)+(d-3)k-2s_0+2s_2\\
&&=n(d-3/2)-2k+2s_1+4s_2.
\end{eqnarray*}}
Hence the number of ways to complete the pairing is
\begin{equation}\label{eq22}
N\left(n(d-3/2)-2k+2s_1+4s_2\right).
\end{equation}

Multiply (\ref{eq17a}) -- (\ref{eq18}) by (\ref{eq21}) --
(\ref{eq22}), sum over all possible $S_0$, then multiply by
(\ref{eq17}). This results in the number of pairings containing a
3-star factor and an oriented and rooted cycle
\begin{eqnarray}\label{eq23}
&&\sum_{S_0}~\frac{n!}{(n-k)!}~(d(d-1))^k~
{{n-k}\choose{n/4-s_0-s_1-s_2}}~(d(d-1)(d-2))^{n/4-s_0-s_1-s_2}\nonumber \\
&&\times~(d-2)^{s_0+s_1+s_2}~(d-3)^{s_0+s_1}(d-4)^{s_0}~\frac{(3n/4-s_1-2s_2)!}{(3!)^{n/4-s_1-s_2}}~(d-2)^{k-s_0-2s_1-3s_2}\nonumber\\
&&\times~d^{3n/4-k+s_0+s_1+s_2}~N\left(n(d-3/2)-2k+2s_1+4s_2\right).
\end{eqnarray}

\noindent Dividing (\ref{eq23}) by the number of pairings with a
3-star factor, which is
\[\frac{n!}{(n/4)!}~\left(d~\left(\frac{(d-1)(d-2)}{3!}\right)^{1/4}\right)^n~N(d(n-3/2))\]
and then evaluating asymptotically, we obtain
\begin{equation}\label{eq24}
\left(\frac{3(d-1)(d-2)}{4(d-3/2)}\right)^k \sum_{S_0} \left(\frac{(d-3)(d-4)}{3(d-1)(d-2)}\right)^{s_0} \left(\frac{8(d-3/2)(d-3)}{3(d-1)(d-2)^2}\right)^{s_1}
\left(\frac{32(d-3/2)^2}{9(d-1)(d-2)^3}\right)^{s_2}.
\end{equation}

We follow an approach used in \cite{Jan95} to determine the summation. We can view \emph{$\textbf{0}$, $\textbf{1}$,
$\textbf{2}$}~~as three states in a Markov Chain, where the final state is equal to the initial state. We observe that
\begin{itemize}
\item \emph{$\textbf{1}$} followed by \emph{$\textbf{0}$} means we
pass a 1-star and this contributes a factor
\[\frac{8(d-3/2)(d-3)}{3(d-1)(d-2)^2},\]
\item \emph{$\textbf{1}$} followed by \emph{$\textbf{2}$} means we
pass a 2-star and this contributes a factor
\[\frac{32(d-3/2)^2}{9(d-1)(d-2)^3}.\]
\end{itemize}
Thus for the transition matrix given by
\[\mathbf{A}=\left( \begin{array}{ccc}
1+\displaystyle{\frac{(d-3)(d-4)}{3(d-1)(d-2)}}&1&0\\
\displaystyle{\frac{8(d-3/2)(d-3)}{3(d-1)(d-2)^2}}&0&\displaystyle{\frac{32(d-3/2)^2}{9(d-1)(d-2)^3}}\\
1&0&0\\
\end{array}\right)\]
we have
\[Tr(\mathbf{A}^k)= \sum_{S_0} \left(\frac{(d-3)(d-4)}{3(d-1)(d-2)}\right)^{s_0}
 \left(\frac{8(d-3/2)(d-3)}{3(d-1)(d-2)^2}\right)^{s_1}\left(\frac{32(d-3/2)^2}{9(d-1)(d-2)^3}\right)^{s_2}.\]
Since the eigenvalues of $\mathbf{A}$ are
\begin{eqnarray*}
\gamma1&=&\frac{4(d-3/2)}{3(d-2)},\\
\gamma2&=&\frac{-3(d-2)+\sqrt{-15d^2+24d}}{3(d-1)(d-2)},\\
\gamma3&=&\frac{-3(d-2)-\sqrt{-15d^2+24d}}{3(d-1)(d-2)},
\end{eqnarray*}
then
\[Tr(\mathbf{A}^k)=\left(\frac{4(d-3/2)}{3(d-2)}\right)^k + \left(\frac{-3(d-2)+\sqrt{-15d^2+24d}}{3(d-1)(d-2)}\right)^k + \left(\frac{-3(d-2)-\sqrt{-15d^2+24d}}{3(d-1)(d-2)}\right)^k.\]
Then (\ref{eq24}) becomes
\begin{equation}\label{eq25}
(d-1)^k~+~\left(\frac{-3(d-2)+\sqrt{-15d^2+24d}}{4(d-3/2)}\right)^k~+~\left(\frac{-3(d-2)-\sqrt{-15d^2+24d}}{4(d-3/2)}\right)^k.
\end{equation}
After (\ref{eq25}) is divided by $2k$ to remove the orientation and rooting of the cycle we obtain (\ref{eq16}). $\square$

\section{Proof of Theorem \ref{main}}
First we prove the following theorem.
\begin{t1}\label{theo4}
Let $n\equiv 0$ mod $4$ and $4 \leq d \leq 10$. Then for $P\in
\mathcal{P}_{n,d}$, $G(P)$ a.a.s has a 3-star factor. Moreover,
\[\frac{Y_d^*}{\textbf{E}Y_d^*}\rightarrow W =\prod_{k=1}^\infty(1+\delta_k)^{Z_k} e^{-\lambda_k \delta_k}\mbox{ for } n \rightarrow \infty\]
where $Z_k$ are independent Poisson variables with $\textbf{E}Z_k = \lambda_k$ for $k\geq 1$ and
\begin{eqnarray*}
\lambda_k&=&\frac{(d-1)^k}{2k},\\
\delta_k&=&\left(\frac{-3(d-2)+\sqrt{-15d^2+24d}}{4(d-1)(d-3/2)}\right)^k
+\left(\frac{-3(d-2)-\sqrt{-15d^2+24d}}{4(d-1)(d-3/2)}\right)^k.
\end{eqnarray*}
\end{t1}
\textbf{Proof.} We will show that $Y_d^*$ satisfies the conditions (A.1) -- (A.4) in (Theorem 4.1, \cite{Wor99}).
Since $X_k$ is the number of short cycles of length $k$ in a pseudograph coming from ${\cal P}_{n,d}$,
then (A.1) is fulfilled with $\lambda_k=\frac{(d-1)^k}{2k}$,
by Bollob\'as' result on short cycles in ${\cal P}_{n,d}$ \cite{Bol85}.
The condition (A.2), (A.3) and (A.4) are fulfilled consecutively by Lemma \ref{le4},
Theorem \ref{expectation} and Theorem \ref{variance}. $\square$ \\

\noindent \textbf{Proof of Theorem \ref{main}.}\\
Theorem \ref{main} comes directly from Theorem \ref{theo4} by
Lemma \ref{pairing}. From the argument in (Remark 9.25,
\cite{Jan00}) we also obtain
\begin{eqnarray*}
\frac{\textbf{E}Y_d}{\textbf{E}Y_d^*}&\rightarrow&\exp\left(\frac{3(5d^2-12d+6)}{4(2d-3)^2}\right), \\
&&\\
\frac{\textbf{E}Y_d^2}{(\textbf{E}Y_d)^2}&\rightarrow&\exp\left(\frac{-9(8d^5-63d^4+206d^3-322d^2+216d-36)}{4(2d-3)^4(d-1)^2}\right)\\
&&\times~\left(\frac{2(d-1)^{1/2}(d-3/2)^2}{(d-3)(4d^3-13d^2+36d-36)^{1/2}}\right).~~~~~~~~~~~~~~~~~~~~~~~~~~~~~~~~~~\square
\end{eqnarray*}

We should point it out again, that most part of the computation in
determining the second moment of the number of 3-star factors in
$G\in \mathcal{P}_{n, d}$ is valid for general $d \geq 4$. The only
part that is still hard to prove is showing that the desired maximum
point is the global maximum. From what we have in the case $4 \leq d
\leq 10,$ we conjecture that the asymptotic distribution of the
number of 3-star factor in random $d$-regular graph have the same
behaviour for $d \geq 4$.

\begin{conj}
Restrict $n$ to $0$ mod $4$ and $d \geq 4$. Then $G\in
\mathcal{G}_{n,d}$ a.a.s has a 3-star factor. Furthermore, letting
$Y_d$ denote the number of 3-star factors in $G \in
\mathcal{G}_{n,d}$,
\[\frac{Y_d}{\textbf{E}Y_d}\rightarrow W =\prod_{k=3}^\infty(1+\delta_k)^{Z_k}
 e^{-\lambda_k \delta_k}\mbox{ for } n \rightarrow \infty\]
where $Z_k$ are independent Poisson variables with $\textbf{E}Z_k = \lambda_k$ for $k\geq 3$ and
\begin{eqnarray*}
\lambda_k & = & \frac{(d-1)^k}{2k},\\
\delta_k & = & \left(\frac{-3(d-2)+\sqrt{-15d^2+24d}}{4(d-1)(d-3/2)}\right)^k
+\left(\frac{-3(d-2)-\sqrt{-15d^2+24d}}{4(d-1)(d-3/2)}\right)^k.
\end{eqnarray*}
\end{conj}

\end{document}